\documentclass[11pt]{article}

\usepackage{amsmath}
\usepackage{amssymb}

\newtheorem{thm}{Theorem}[section]
\newtheorem{prop}[thm]{Proposition}

\newtheorem{rem}[thm]{Remark}
\newtheorem{frem}[thm]{Final Remark}

\def\lsp{\Gamma\backslash X}

\def\P{\mathbf P}
\def\ga{\Gamma}
\def\G{\mathbf G}
\def\Z{\mathbb Z}
\def\Q{\mathbb Q}
\def\R{\mathbb R}

\title{Infinite topology of curve complexes and non-Poincar\'e duality
of Teichm\"uller modular groups}
\author{Nikolai Ivanov and Lizhen Ji}
\date{}
\begin{document}

\maketitle

\section{Introduction}

Let $S$ be an orientable surface.
Let $\text{Diff}^{\pm}(S)$ be the group of all diffeomorphisms of $S$,
and $\text{Diff}^{0}(S)$ its identity component.
Then $Mod^{\pm}_ {S}=\text{Diff}(S)/\text{Diff}^{0}(S)$
is called the {\em extended mapping class group}
or the {\em extended Teichm\"uller modular group} of $S$. 
Let $\text{Diff}(S)$ be the subgroup
of orientation preserving diffeomorphisms of $S$.
Then $Mod(S)=\text{Diff}(S)/\text{Diff}^{0}(S)$  
is called the {\em mapping class group} or the {\em Teichm\"uller modular group} of $S$.

When $S$ is a closed surface
of genus 1, for example, when $S=\Z^{2}\backslash \R^{2}$ is the standard torus, then $Mod^{\pm}(S)$ can be identified
with $GL(2, \Z)$, and
$Mod(S)$ can be identified with $SL(2, \Z)$.
If $S$ is a closed surface of genus $g\geq 2$, or more generally
an oriented surface of negative Euler characteristic $\chi(S)$, then
$Mod^{\pm}(S)$ and $Mod(S)$ can be considered as natural generalizations
of $GL(2, \Z)$ and $SL(2, \Z)$ respectively.
 
The group $SL(2, \Z)$ is an important example
of arithmetic subgroups of semisimple linear algebraic groups, in particular, it is the first group in the classical family
of arithmetic subgroups $SL(n, \Z)$, $n\geq 2$ of $SL(n, \R)$.
Arithmetic subgroups $\ga$ of semisimple linear algebraic groups $\G$ (defined over
$\Q$) enjoy many good properties and are special among
discrete groups. For example, they are finitely presented and enjoy
other finiteness properties such as being of type $FP_{\infty}$ and of type $FL$. 
Borel and Serre showed in \cite{bos} that arithmetic subgroups $\ga$
are virtual duality groups and their
virtual cohomological dimension can be computed explicitly (we outline parts of this theory below).
The cohomology groups of a natural family of arithmetic subgroups such as $SL(n, \Z)$
stabilize as $n\rightarrow \infty$.
See \cite{se} and \cite{bo1} for a summary and \cite{bo2} for a computation of the stable real cohomology groups. 
  
Let $G=\G(\R)$ be the real locus of $\G$, let $K\subset G$ be a maximal compact subgroup,
and let $X=G/K$ be the associated symmetric space. Assume that $\ga$ is a torsion-free arithmetic
subgroup of $\G(\Q)$. Then the locally symmetric space $\lsp$ is an aspherical
manifold with the fundamental group $\ga$. 
If $\ga$ is a  cocompact subgroup of $G$, then the Poincar\'e duality for $\lsp$
implies that $\ga$ is a Poincar\'e duality group.
On the other hand, if 
the arithmetic subgroup $\ga$  is not a 
cocompact discrete subgroup of $G$ (as is $SL(n, \Z)$, for example), then \cite{bos} (see the notes by Serre near the end 
\cite{bie}) implies that
$\ga$ is {\em not a Poincar\'e duality group}, in particular, it can not be realized
as the fundamental group of a closed aspherical manifold. Instead of being a Poincar\'e duality group, $\ga$ enjoys a weaker property being 
a {\em duality group} in the sense of Bieri-Eckmann \cite{bie}.
In Section 2,  we will recall the definition of duality groups
and Poincar\'e duality groups (see, in particular, the formulas (\ref{duality-definition}), (\ref{Poincare-duality}), \ref{virtual-duality}\.) .

In the context of Teichm\"uller modular groups, the {\em Teichm\"uller space} $T_S$ of $S$ plays the role of the symmetric space $X$,
and the canonical action of $Mod_S$ on $T_S$ plays the role of the action of $\Gamma$ on $X$. This analogy was discovered by Harvey \cite{harv1}, \cite{harv} who, in particular, was motivated by
the problem of providing analogues of some constructions of Borel-Serre \cite{bos}.
 
In a series of remarkable papers \cite{ha2}--\cite{ha4} (see also \cite{iv5} for another approach to these and other related results, and \cite{ha1}, \cite{iv1} for expository accounts), Harer, motivated by the analogy between the groups $Mod_{S}$ and the arithmetic groups (which was also originally pointed out by Harvey \cite{harv1}, \cite{harv}) established many properties of $Mod_{S}$ similar to those of arithmetic subgroups. 
 
For simplicity, in the following we assume that $S$ is a closed orientable surface of genus $g\geq 2$, and $Mod_{S}$ will be also denoted by $Mod_{g}$. All results are true for surfaces with non-empty boundary also, with minor modifications.
 
We are especially interested in the following result (see \cite{ha2}, Theorem 4.1, and \cite{iv5}, Theorem 6.6 for a different proof).
 
\begin{thm}\label{duality-group}
For every torsion-free subgroup $\ga\subset Mod_{g}$ of finite index,
$\ga$ is a duality group (in the sense of Bieri-Eckmann {\em \cite{bie}}) of dimension $4g-5$. Therefore, $Mod_{g}$
is a virtual duality group of virtual cohomological  dimension $4g-5$.
\end{thm}
  
This theorem is an analogue of some results of Borel-Serre \cite{bos}. 
The crucial ingredient of the proof of these results of Borel-Serre is the Solomon-Tits theorem about the homotopy type
of the spherical Tits building $\Delta(\G)$, which is (as proved by Borel-Serre) is homotopy equivalent
to the boundary of the partial Borel-Serre compactification
of $X$ constructed in \cite{bos}. An analogue of the Borel-Serre compactification for the Teichm\"uller spaces was suggested by Harvey \cite{harv},
as also an analogue of the Tits building $\Delta(\G)$. The latter analogue is the {\em complex of curves} $\mathcal C(S)$ of $S$, which plays a
fundamental role in the theory of Teichm\"uller modular groups $Mod_S$ and, in particular, in the present paper. See the beginning of
Section 2 for the definition of $\mathcal C(S)$.

The underlying topological space of the spherical Tits building $\Delta(\G)$ 
consists of {\em non-disjoint union} of countably infinitely many spheres,
and hence it is not surprising that $\Delta(\G)$ is homotopy equivalent to 
a bouquet of countably infinitely many spheres, as it is the case by the Solomon-Tits Theorem.
Basically, this is proved by showing the intersection of any two distinct
apartments is contractible (see \cite{br1}, p. 92, Theorem 2, and \cite{ku}).
On the other hand, the curve complex $\mathcal C(S)$ does not have any structure
of apartments, and its homotopy type is not clear (or easy to be guessed)
from the definition. 

In  \cite{ha2} (see also \cite{ha1}, Chap. 4, \S 1, and \cite{iv5}, \S 3 or \cite{iv1}, Theorem 3.3.A for a different proof), Harer
proved the following analogue of the Solomon-Tits theorem.

\begin{thm}\label{bouquet}
The curve complex $\mathcal C(S)$ is homotopy equivalent to a bouquet of spheres
$\vee S^{n}$, where the dimension $n=2g-2$. 
\end{thm}

The proof of the above theorem \ref{bouquet} given in \cite{ha2}, \cite{ha1} 
goes in two steps: (1) $\mathcal C(S)$ is $(n-1)$-connected;
(2) the homotopy dimension of $\mathcal C(S)$ is bounded from above by $n$.
Therefore, $\mathcal C(S)$ is homotopy equivalent to a  bouquet $\vee S^{n}$.
The alternative proof of \cite{iv1}, \S 3.3, p. 546, \cite{iv5}, \S 3 follows the same general outline.
These arguments leave open the possibility that the number of spheres in the above bouquet is equal to $0$ and and $\mathcal C(S)$
is contractible. In fact, this is not the case, as the following theorem shows.

\begin{thm}\label{curve-nontrivial}
The curve complex $\mathcal C(S)$ is not contractible, and hence
the bouquet of spheres $\vee S^{n}$ in Theorem \ref{bouquet}
contains at least one sphere.
\end{thm}

Motivated by Theorems \ref{bouquet}, \ref{curve-nontrivial}, and by the Solomon-Tits Theorem, it is also natural expect
the following improvement of Theorem \ref{bouquet}.

\begin{thm}\label{infinite-topology}
The curve complex $\mathcal C(S)$ is homotopy equivalent to the bouquet
of countably infinitely many spheres $S^n$, where $n=2g-2$.
\end{thm}

The algebraic counterpart of Theorems \ref{curve-nontrivial}, \ref{infinite-topology} is the following result.

\begin{thm}\label{non-poincare}
For every torsion-free subgroup $\ga\subset Mod_{S}$ of finite index,
$\ga$ is not a Poincar\'e duality group.
Equivalently, $Mod_{S}$ is not a virtual Poincar\'e duality group.
In particular, no subgroup $\ga$ of finite index of $Mod_S$ can be realized as
the fundamental group of a closed aspherical manifold.
\end{thm}

If we take seriously the analogy between the modular groups $Mod_S$ and the arithmetic groups, and, in particular, the analogy
between the Techm\"uller space $T_S$ of $S$ and the symmetric spaces $X=G/K$, this result is to be expected. Namely, the arithmetic groups $\Gamma$ are virtual Poincar\'e duality groups if and only if they act on $X$ cocompactly, 
and the action of $Mod_S$ on $T_S$ is not, as is well known.
This heuristic argument is based only on the fact the the complex of curves $\mathcal C(S)$ is non-empty (this is equivalent to the action on $T_S$ being non-cocompact). 
But the actual proof proceeds differently. Namely, the proof of Theorem \ref{infinite-topology} is based on Theorem \ref{non-poincare}, 
the proof of which is, in turn, based the properties of the Mess subgroups \cite{me} of $Mod_S$. The proof of Theorem \ref{curve-nontrivial} is more
direct, but it is also uses properties of the group $Mod_S$; namely, the fact that $Mod_S$ has non-zero (virtual) cohomological dimension.

Theorem  \ref{curve-nontrivial}, while undoubtedly known to the workers in the field 
(including J. Harer and the first author) for a long time, was not stated or proved explicitly in the literature.
At the same time it was used in several places, for example, in understanding relations 
between the automorphism group of the curve complex $\mathcal C(S)$ and the mapping class group of $S$
(see \cite{iv4}, \cite{kor}, \cite{kor2}, and \cite{luo}), and was also used implicitly in the proof of Theorem \ref{duality-group}.
We included the relevant part of the proof of Theorem \ref{duality-group} at the end of Section 3.

While Theorem \ref{curve-nontrivial} is at least implicitly contained in the papers \cite{ha2}, \cite{iv5},
Theorems \ref{infinite-topology} and \ref{non-poincare} are not discussed there, and the proof of these theorems which we offer
below is based on later results due to Mess \cite{me}. One can safely say that everybody assumed that these results are
true, but nobody cared to write down a proof. (After a preliminary version of this note was written up, T. Farrell pointed out that
Harer claimed in \cite{ha1}, p. 180, lines 20-22, that the dualizing module $I$ (denoted by $C$
in this note) for a torsion-free subgroup $\ga$ of $Mod_S$ of finite index has infinite rank,
which is equivalent to Theorem \ref{infinite-topology} of this note.)

This question on the nontrivial topology of $\mathcal C(S)$ has been raised by several people
and has caused some confusions. The goal of this note is provide proofs of Theorems \ref{curve-nontrivial}, 
\ref{infinite-topology}, and \ref{non-poincare}. 

The proof is based on a consideration of
cohomology groups of torsion-free subgroups of $Mod_S$ of finite index
and on the relation the curve complex $\mathcal C(S)$ to the boundary structure
of a Borel-Serre type partial compactification $\overline{T_{S}}^{BS}$ discovered by Harvey \cite{harv} (or a truncation $T_{S}(\varepsilon)$) 
of the Teichm\"uller space $T_S$ of $S$. These ideas form the foundation of the theory presented in \cite{ha2}, \cite{ha1}, and \cite{iv5}.
In particular, to a big extent our proof of Theorem \ref{curve-nontrivial} is
already contained in these papers. 

In order to have a fuller analogy between the arithmetic and the Teichm\"uller modular cases, 
it seems desirable to have a proof of Theorem \ref{curve-nontrivial} using internal structure of $\mathcal C(S)$
as in the case of the Solomon-Tits  theorem for spherical Tits buildings.
Note that cohomology of arithmetic groups and the action of arithmetic
subgroups on the spherical Tits buildings $\Delta(\G)$ are not used to
prove the Solomon-Tits Theorem. So the proof of Theorem \ref{curve-nontrivial} presented below looks like putting
a cart in front of the horse. It is also worthwhile to understand if there are some special subcomplexes
of $\mathcal C(S)$ which play the role of apartments in Tits buildings.\\

\noindent{\em Acknowledgments.} This work started during a workshop
{\em Aspects of Teichm\"uller spaces} at the University of Michigan in April 2007
when Feng Luo brought the lack
of a proof of Theorem \ref{curve-nontrivial} to the attention of the second author.
This problem was also pointed out by Yair Minksy and Benson Farb. The workshop
{\em Aspects of Teichm\"uller spaces} at University of Michigan (April 12-14, 2007)
was supported by the RTG grant  0602191 in {\em Geometry, Topology and Dynamics}, 
and the FRG grant {\em Geometric Function Theory}.

The first author would like to thank the organizers of this workshop (and the second
author among them) for the invitation to speak at the workshop and for the stimulating
atmosphere.

The second author would like to thank Feng Luo 
also for comments on an earlier version;
he would like to thank Tom Farrell for conversations and references on duality groups,  
and also for other valuable advice on writing; 
he would also like to thank S. Wolpert and R. Spatzier for helpful
conversations and comments, and J. Harer, B. Farb and K. Brown for  helpful email correspondences.
 
The first author was partially supported by the NSF Grant DMS-0406946, 
and the second author was partially supported by NSF Grant DMS 0604878.

\section{Preliminaries} 

For simplicity, we assume that $S$ is a closed oriented surface of genus $g\geq 2$.
By the definition, the vertices of the curve complex $\mathcal C(S)$
are free homotopy classes $\langle c \rangle$ of simple closed curves $c$ in $S$.
Vertices $\langle c_{1}\rangle, \ldots, \langle c_{k+1}\rangle$
form the vertices of a $k$-simplex if and only if they are all different and every two curves $c_{i_{1}}$
and $c_{i_{2}}$ for $1\leq i_1<i_2\leq k+1$ are isotopic to disjoint curves.
It is well known and easy to see that $\mathcal C(S)$ is a simplicial complex of dimension $3g-4$. 

Let $T_S$, which is also denoted by $T_g$, be the Teichm\"uller space of hyperbolic metrics
on $S$. The curve complex $\mathcal C(S)$ was first introduced by Harvey \cite{harv} in order to
construct a partial compactification $\overline{T_g}^{BS}$ of $T_g$ (his project was completed by the first author in \cite{iv2}), 
which is similar to and motivated by the Borel-Serre
partial compactification $\overline{X}^{BS}$ of symmetric spaces $X=G/K$.
The latter is a manifold with corners whose boundary faces (or components) are parametrized
by simplexes of a spherical Tits building associated with $G$. 

More specifically, if $G=\G(\R)$ is the real locus of a linear semisimple algebraic group
$\G$ defined over $\Q$, then the boundary of $\overline{X}^{BS}$ is decomposed into 
boundary components $e(\P)$
parametrized by $\Q$-parabolic subgroups $\P$ of $\G$ satisfying the conditions:
\begin{enumerate}
\item The closure of every $e(\P)$ in $\overline{X}^{BS}$ is contractible.
\item $e({\P_1})$ is contained in the closure of $e({\P_2})$ if and only if 
$\P_1$ is contained in $\P_2$.
\end{enumerate}

Recall that the spherical Tits building $\Delta(\G)$ of the algebraic group $\G$ defined over $\Q$
is a simplicial complex with one simplex $\sigma_\P$ for each proper parabolic
subgroup $\P$ (defined over $\Q$) of $\G$ satisfying the conditions:
\begin{enumerate}
\item Vertices of $\Delta(\G)$ correspond to maximal proper $\Q$-parabolic subgroups of $\G$.
\item $\sigma_{\P_1}$ is contained in $\sigma_{\P_2}$ if and only if $\P_2$ is contained in
$\P_1$. In particular maximal proper $\Q$-parabolic subgroups $\P_1, \cdots, \P_{k+1}$
form the vertices of a $k$-simplex if and only if their intersection $\P_1\cap \cdots \cap \P_{k+1}$
is a parabolic subgroup of $\G$. 
\end{enumerate}

The boundary components $e(\P)$ of $\overline{X}^{BS}$ correspond to the simplexes
of the Tits building $\Delta(\G)$ of $\G$. Furthermore, 
the partial Borel-Serre compactification $\overline{X}^{BS}$ is a real analytic manifold with corners
whose boundary $\partial \overline{X}^{BS}$ has the same homotopy type as the (geometric realization of the) building $\Delta(\G)$.

For every torsion-free arithmetic subgroup $\ga\subset \G(\Q)$, the action of $\ga$ on $X$
extends to a proper and continuous action on $\overline{X}^{BS}$ such that the quotient
$\ga\backslash \overline{X}^{BS}$ is a compact manifold with corners, in particular, is a finite
$B\ga$-space, i.e., a classifying space $B\ga$ which is a finite CW-complex. 
These results are used effectively in studying cohomology groups of $\ga$.
In fact, it is shown in \cite{bos} (see also \cite{bie}) that 
every arithmetic subgroup $\ga$ of $\G(\Q)$
is a virtual duality group with virtual cohomological dimension equal to $\dim X-r$,
where $r$ is equal to the $\Q$-rank of $\G$. This result motivated the work \cite{ha2} and other related
results on mapping class groups (see \cite{ha1} for a summary).

Note that the dimension of $\Delta(\G)$
is equal to $r-1$.  It is a well-known theorem of Solomon-Tits (see \cite{br1}, p. 93, Theorem 2) that
$\Delta(\G)$ is homotopy equivalent to a bouquet of countably infinitely many
spheres of dimension ${r-1}$.
In fact, if we fix a chamber $C$ of the building $\Delta(\G)$,
there is one sphere for each apartment of the building $\Delta(\G)$ containing $C$.
This implies that when the $\Q$-rank $r$ is positive, every arithmetic subgroup 
of $\G(\Q)$ is not a Poincar\'e duality group,
in particular, can not be realized as the fundamental group of a closed aspherical  
manifold. (See Note 1 by J.P.Serre on  \cite{bie}, p.124.)  Note also that
the condition that the $\Q$-rank $r>0$ is equivalent to that every arithmetic subgroup
$\ga$ of $\G(\Q)$ is not a cocompact  discrete subgroup of $G$. 
 
For the convenience of the reader, we recall (see, for example \cite{bie}) that 
a discrete group $\ga$ is called {\em a Poincar\'e duality group} of dimension $n$
there exists  a right $\ga$-module structure on $\Z$ such that 
\begin{equation}\label{duality-definition}
H^k(\ga, A)\cong H_{n-k}(\ga, \Z\otimes A)
\end{equation}
for all $k$ and all $\ga$-modules $A$.
In this case,
\begin{equation}\label{Poincare-duality}
H^n(\ga, \Z\ga)\cong \Z.
\end{equation}
If  $\ga$ admits a classifying space which is a closed (i.e. compact  without boundary) manifold $M$
or equivalently if $\ga$ is the fundamental
group of a closed aspherical manifold $M$, then the Poincar\'e duality  for $M$
implies that $\ga$ is a Poincar\'e duality group.

A discrete group $\ga$ is called a {\em duality group}  of dimension $n$ (in the sense of Bieri-Eckmann \cite{bie})
with respect to a right $\ga$-module $C$, called the {\em dualizing module}, if there is an element $e\in H_n(\ga, C)$
such that the cap-product with $e$ induces isomorphisms
\begin{equation}\label{virtual-duality}
H^k(\ga, A)\cong H_{n-k}(\ga, C\otimes A)
\end{equation}  
for all $k$ and all left $\ga$-modules $A$.
In this case, the {\em cohomological dimension} of $\ga$ is equal to $n$. 

Let us turn to the Teichm\"uller case. The Borel-Serre
partial compactification $\overline{T_g}^{BS}$ of \cite{harv}
 is also a manifold with corners with  boundary faces
$B_\sigma$ parametrized by simplexes $\sigma$ of the curve complex $\mathcal C(S)$ such that 
$B_{\sigma_1}$ is contained in the closure of $B_{\sigma_2}$ if and only if
$\sigma_2$ is contained in $\sigma_1$. 
In this sense, the curve complex $\mathcal C(S)$ 
is similar to spherical Tits buildings $\Delta(\G)$. 
The construction of  $\overline{T_g}^{BS}$ was not worked
out in detail in \cite{harv} and $\overline{T_g}^{BS}$ was later constructed in \cite{iv2}.

A version of $\overline{T_g}^{BS}$ which is easier to construct and to deal with is provided by a truncation $T_g(\varepsilon)$
by the Teichm\"uller space $T_g$. It is 
due to the first author \cite{iv3} (see also \cite{iv1}, Theorem 5.4.A). 
For every point $x=(S, d)\in T_g$, where $d$ is a hyperbolic metric on $S$,
 denote the length of a closed geodesic $c$
in $S$ with respect to the hyperbolic metric $d$ by $\ell_x(c)$.
For $\varepsilon$ sufficiently small,
define 
$$T_g(\varepsilon)=\{ x=(S, d)\in T_g\mid \text{\ for every closed geodesic\ }
 c\subset S,\quad 
\ell_x(c)\geq \varepsilon \}.$$
The smallness of $\varepsilon$ is required so that two geodesics of lengths less than $\varepsilon$
do not intersect. Such $\varepsilon$ always exists due to the collar theorem
(see \cite{bu}, \S 4.1 for example).

Let $Mod_g=Mod_S$ be the Teichm\"uller modular group of $S$.
Then $Mod_g$ acts properly on $T_g$. Since points of every orbit of $Mod_g$
in $T_g$ represent isometric metrics on $S$, it is clear that $Mod_g$ leaves $T_g(\varepsilon)$
invariant. It follows from the Mumford compactness criterion
that the quotient $Mod_g\backslash T_g(\varepsilon)$ is compact.

The following result (see \cite{iv1}, Theorem 5.4.A) is crucial for us.

\begin{prop}\label{boundary}
The space $T_g(\varepsilon)$ is a contractible manifold with corners.
Its boundary faces when $T_g(\varepsilon)$ is considered 
as a manifold with corners are contractible and parametrized by the simplexes of the curve complex
$\mathcal C(S)$, and the whole boundary $\partial T_g(\varepsilon)$ is homotopy equivalent  
to $\mathcal C(S)$. In particular, for every torsion-free subgroup $\ga\subset Mod_g$,
the quotient $\ga\backslash T_g(\varepsilon)$ is a finite $B\ga$-space. 
\end{prop}

The crucial part of the proof of this theorem is to show that $T_g(\varepsilon)$ is contractible
and this follows from the fact $T_{g}(\varepsilon)$ is a deformation retract of $T_{g}$.
The correspondence between its boundary components
and the simplexes of $\mathcal C(S)$ is described as follows.
Each simplex $\sigma$ of $\mathcal C(S)$ gives a collection of disjoint and non-isotopic simple closed geodesics. 
The corresponding boundary face of $\partial T_g(\varepsilon)$
consists of those marked hyperbolic metrics on $S$ where the geodesics in this collection
have length exactly
equal to $\varepsilon$. Therefore, each boundary face of $\partial T_g(\varepsilon)$
is basically the product of a truncated Teichm\"uller space and an Euclidean space (corresponding to the Fenchel-Nielsen
twist parameters),
and hence also contractible. 

\section{Proofs of Theorems \ref{curve-nontrivial}, \ref{infinite-topology}, and \ref{non-poincare}}

We start with the proof of Theorem \ref{curve-nontrivial}. It is similar to the arguments in \cite{ha2}, \cite{ha1}, \cite{iv5},  \cite{iv1}, 
and, of course, to the arguments in \cite{bos} where Borel and Serre prove that arithmetic subgroups of linear semisimple algebraic groups are
virtual duality groups.\\

\noindent{\em Proof of Theorem \ref{curve-nontrivial}. } Let $\ga\subset Mod_g$ be a torsion-free subgroup of finite index.
The idea of the proof is to show that if $\mathcal C(S)$ is homotopy trivial,
i.e., there is no $S^n$ in the bouquet $\vee S^{n}$ of spheres in Theorem \ref{bouquet},
then the cohomological dimension of $\ga$, denoted by ${\rm cd}\,\ga$, is equal to zero.
But this is impossible since $\ga$ is a torsion-free infinite group.

By Proposition \ref{boundary},  $\ga\backslash T_g(\varepsilon)$ is a compact $B\ga$-space.
Hence by \cite{iv1}, Corollary 6.1.G, 
${\rm cd}\,\ga$ is equal to the maximum number $n$ such that 
$${H}_c^n( T_g(\varepsilon),\Z)\neq 0,$$
where $H_c^n(\,\cdot\,)$ denotes the cohomology with compact supports. 
Denote the dimension of $T_g(\varepsilon)$ by $d$.
Following \cite{iv1}, Theorem 6.1.H (or \cite{bos}),
we see that by the Poincar\'e-Lefschetz duality,
\begin{equation}\label{proof-1}
{H}_c^n(T_g(\varepsilon), \Z)=H_{d-n}(T_g(\varepsilon),  \partial T_g(\varepsilon); \Z).
\end{equation}
Since $T_g(\varepsilon)$ is contractible, the homology long exact sequence gives
\begin{equation}\label{proof-2}
H_{d-n}(T_g(\varepsilon),  \partial T_g(\varepsilon); \Z)
\cong \widetilde{H}_{d-n-1}(\partial T_g(\varepsilon), \Z).
\end{equation}
By Proposition \ref{boundary},
$\partial T_g(\varepsilon)$ is homotopy equivalent to $\mathcal C(S)$.
Therefore,
\begin{equation}\label{eq-curve-complex}
{H}_c^n(T_g(\varepsilon), \Z)
\cong \widetilde{H}_{d-n-1}(\mathcal C(S)).
\end{equation}
Suppose that $\mathcal C(S)$ is homotopy equivalent to a point. Then 
Equation (\ref{eq-curve-complex}) implies that for every  $n\geq 1$,
$$ {H}_c^n(T_g(\varepsilon), \Z)=0.$$
This in turn implies that ${\rm cd}\,\ga=0$.

It is well known that for every pairs of groups $\ga'\subset \ga$,  ${\rm cd}\,\ga'\leq{\rm cd}\,\ga$.
(See \cite{br1}, Proposition 2.4, p. 187 or  \cite{iv1},  \S 6.4, p. 584.)
Since $\ga$ is torsion-free and infinite, it contains $\Z$ as a subgroup.
A classifying space of $\Z$ is given by the circle $S^1$, and hence clearly  ${\rm cd}\,\Z=1$.
This implies that ${\rm cd}\,\ga\geq 1$.
This contradiction implies that the homotopy type of $\mathcal C(S)$ cannot be trivial.
The proof of  Theorem  \ref{curve-nontrivial} is complete. $\Box$

\begin{rem}\label{one} {\em
The proof of Theorem \ref{curve-nontrivial} can be arranged in slightly different way.
By Theorem \ref{bouquet} and the arguments based on equations (\ref{proof-1}),  
(\ref{proof-2}) and (\ref{eq-curve-complex}), we see that for $i\neq \nu=4g-5$, $H^i(\ga, \Z\ga)=0$.
If $H^\nu(\ga, \Z\ga)=0$, then \cite{br1}, Proposition 6.7 (see p. 202) implies that
the cohomological dimension ${\rm cd}\,\ga$  is equal to 0.

Since $\ga$ is torsion-free and contains $\Z$, the last paragraph of
the above proof  \ref{curve-nontrivial} implies ${\rm cd}\,\ga\geq 1$.
This, in turn, implies that $H^\nu(\ga, \Z\ga)\neq 0$.
By the above equations again, it follows that the bouquet $\vee S^n$
for $\mathcal C(S)$ contains at least one sphere.

Or, alternatively, 
by \cite{br1}, Theorem 10.1 (see p. 220), and its proof, this non-vanishing
result implies that $\ga$ is a duality group of dimension $\nu$.
Since the lower dimensional cohomology groups of $\ga$ is known to be non-zero (see 
\cite{ha4} for example), the dualizing module $\widetilde{H}_n(\mathcal C(S), \Z)$
is not 0, which implies that $\mathcal C(S)$ is not contractible.}

The last approach has the disadvantage of being dependent on very difficult calculations of Harer \cite{ha4}.

\end{rem}

\begin{rem} \label{FP} {\em
We also note that  although it is not stated explicitly  in  (c)  Theorem 10.1 of \cite{br2}, 
under the assumption that $\ga$ is a torsion-free infinite group of type $FP$, the non-vanishing of 
$H^{\nu}(\Gamma, \Z\Gamma)$ is implied by other conditions: $H^i(\ga, \Z \ga)=0$
for $i\neq \nu$. In fact, since $\ga$ contains $\Z$,
it follows from the second paragraph of  the previous remark
that  the group $H^{\nu}(\Gamma, \Z\Gamma)$ must be non-zero.
 
This non-vanishing of $H^{\nu}(\Gamma, \Z\Gamma)$ 
is used at several places in the book of Brown \cite{br2}. For example, to show that the cohomological dimension of $\Gamma$ is equal to $n$
in \cite{br2}, Chapter VIII, Proposition 6.7, in the proof (ii) $\Rightarrow$ (iv) in  \cite{br2}, Chapter VIII, Theorem 10.1,  
and also in \cite{br2}, Exercise 5 (b) in Section VIII.6. 
 
Note that if $\ga$ is a torsion-free subgroup of finite index in $Mod_g$, then the condition that $\ga$ is of type $FP$ follows from the fact that $\ga$ acts freely
on the manifold with corners $T_g(\varepsilon)$ with a compact quotient 
$\ga\backslash T_g(\varepsilon)$
together with the fact that a compact manifold with corners admits a finite triangulation,
which implies that $\ga$ admits a finite $B\ga$-space and hence is of type $FL$, which is stronger
than being of type $FP$.}
\end{rem}

\noindent{\em Proof of Theorem \ref{non-poincare}.} Suppose that $\ga$ is a Poincar\'e duality group.
Then, by a theorem of R. Strebel (see \cite{st}, 1.2), for every subgroup $\ga'\subset \ga$
of infinite index,  ${\rm cd}\, \ga <  {\rm cd}\,\ga$.

Therefore, in order to prove Theorem \ref{non-poincare}, it suffices to find a subgroup
$\ga'$ of $\ga$ of infinite index (i.e. $[\ga : \ga']=+\infty$)
with ${\rm cd}\,\ga'={\rm cd}\,\ga$. 
One can take a Mess subgroup $B_{g}$ \cite{me} as such a group. See \cite{iv1}, \S 6.3 for an accessible exposition of the relevant part of the Mess work.
The Mess subgroups $B_{g}$ are constructed by an induction on $g$. The construction involves a lot of arbitrary choices, so
actually there are (infinitely) many Mess subgroups of $Mod_g$ for any given $g$. The last step of the construction
of $B_g$ is to take a subgroup of $Mod_g$ supported on a subsurface of genus $g-1$ with $1$ boundary component (it is
constructed from a $B_{g-1}$) and to add to it a Dehn twist $t_c$ about a nontrivial circle $c$ disjoint from this subsurface. Clearly,
the subgroup $B_g$ is contained in the centralizer of $t_c$. In particular, it does not contains nontrivial powers of Dehn twists
about circles not isotopic to a circle disjoint from $c$.
It follows that $B_{g}$ is of infinite index in $Mod_g$.
The intersection $\ga'=B_{g}\cap \ga$ is of finite index in $\ga$.
A crucial result about the Mess subgroups $B_{g}$ is that
$B_{g}$ is the fundamental group of a closed topological manifold of dimension equal to the ${\rm cd}\, \ga=4g-5$; see 
\cite{iv1}, \S 6.3, Theorem 6.3.A. In particular,
${\rm cd}\, \ga'={\rm cd}\,B_{g}={\rm cd}\, \ga$. 
This completes the proof  that $\ga$ is not a Poincar\'e duality group. $\Box$\\

\noindent{\em Proof of Theorem \ref{infinite-topology}.} Suppose that $\mathcal C(S)$ is homotopy equivalent to a bouquet of
finitely many spheres $S^n$. 

Let $\nu= 4g-5$ be the virtual cohomological dimension of $Mod_g$. If $d=\dim T_g(\varepsilon)=6g-6$, then $d-\nu-1=
6g-6-(4g-5)-1=2g-2=n$. For any torsion-free subgroup
$\ga$ of $Mod_g$ of finite index, we have by \cite{iv1}, Lemma 6.1.F
$$H^{\nu}(\ga, \Z \ga)\cong H^{\nu}_c(T_{g}(\varepsilon), \Z).$$ 
Further, by the equations (\ref{proof-1}), (\ref{proof-2}) and (\ref{eq-curve-complex}) above (and the fact that $d-\nu-1=n$), we have
\begin{equation}\label{identity}
H^{\nu}(\ga, \Z \ga)\cong H^{\nu}_c(T_{g}(\varepsilon), \Z)\cong 
\widetilde{H}_{n}(\partial T_{g}(\varepsilon), \Z )\cong
\widetilde{H}_{n}(\mathcal C(S), \Z).
\end{equation}
is a finitely generated abelian group. 

By a theorem of Farrell (see \cite{fa} Theorem 3), if $\ga$ is a finitely presented group of type (FP), 
is a duality group of dimension $\nu$, and if $H^{\nu}(\ga, \Z \ga)$ is finitely generated,
then $H^{\nu}(\ga, \Z \ga)$ is cyclic and hence $\ga$ is a Poincar\'e duality group.  
By Theorem \ref{duality-group},
a torsion-free subgroup $\ga$ of $Mod_S$
of finite index is a duality group. By the Remark \ref{FP},
such a subgroup $\ga$ is of type $FP$. 
This implies that $\ga$ is a Poincar\'e duality group.
But this contradicts Theorem \ref{non-poincare}.

This implies that $\mathcal C(S)$ is homotopy equivalent to a bouquet of
infinitely many spheres $S^n$. Since $\mathcal C(S)$ is a countable simplicial
complex, there are only countably infinitely many spheres
in the bouquet $\vee S^n$. $\Box$

\begin{rem}\label{two}{\em
One can avoid using \cite{fa}, Theorem 3 in the proof of Theorem \ref{infinite-topology}
by applying \cite{br1}, Exercise 4 in Section VIII.6. According to this Exercise, if $H^\nu(\ga, \Z\ga)$ is finitely generated (where $\nu={\rm cd}\, \ga$), then 
for every subgroup $\ga'\subset \ga$ of infinite index,  ${\rm cd}\, \ga <  {\rm cd}\,\ga$. Notice that, if there are only finitely many spheres
in the bouquet, then $H^\nu(\ga, \Z\ga)$ is finitely generated. Using the Mess subgroups as in the proof of Theorem
\ref{infinite-topology}, one can complete the proof.}
\end{rem}

\begin{rem}{\em
In order  to prove Theorem \ref{infinite-topology}, i.e. that
the bouquet $\vee S^n$ contains infinitely many spheres,
it would be sufficient to show the bouquet $\vee S^n$ contains at least two spheres.
In fact, by equation (\ref{identity}), in this case
$H^{\nu}(\ga, \Z \ga)$ contains a subgroup isomorphic to $\Z^{2}$.
Since for the Poinca\'re duality groups the equation (\ref{Poincare-duality}) holds, 
this implies that $\ga$ is not a Poincar\'e duality group, and then one would be able to complete the proof by appealing
to the Farrel's theorem again. Unfortunately, it seems that there is no direct proof of the fact that the bouquet $\vee S^n$ contains at least two spheres.
Compare Remark \ref{solomon-tits}.}
\end{rem}

\noindent{\em Sketch of a proof of Theorem \ref{duality-group}.} For the sake of completeness and clarification of any possible confusion,
we outline a proof of Theorem \ref{duality-group}. 
The arguments are the same as \cite{ha2}, \cite{ha1}, \cite{iv5}, and \cite{bos} (see \cite{bos}, Theorem 8.6.5), 
and prove that $Mod_g$ is a virtual duality group of dimension $\nu=(6g-6)-(2g-2)-1=4g-5$.

Specifically, let $\ga$ be a torsion-free subgroup of $Mod_S$ of finite index.
By \cite{bie}, Theorem 2.3 (see also \cite{br2}, Theorem 10.1, p. 220),
it suffices to prove the following claim: 
\begin{equation}\label{claim}
H^n(\ga, \Z\ga)=0  \text{ \  for\  } n\neq \nu, \text{\ and \ } H^{\nu}(\ga, \Z\ga)
\text{\ is {\em nonzero} and torsion-free}.
\end{equation}
Since $\ga\backslash T_g(\varepsilon)$ is a finite $B\ga$-space, by \cite{iv1}, Lemma 6.1.F, we have
\begin{equation}\label{iden}
H^n(\ga, \Z\ga)=H^n(\ga\backslash T_g(\varepsilon), \Z\ga)\cong H^n_c(T_g(\varepsilon), \Z).
\end{equation}
Since the boundary  $\partial T_g(\varepsilon)$ is homotopy equivalent to
a nonempty bouquet of spheres $\vee S^n$, where $n=2g-2$,
the same arguments as in the proof of  Theorem \ref{curve-nontrivial}
above (see equations (\ref{proof-1}),  (\ref{proof-2}) and (\ref{eq-curve-complex})\;)
prove (\ref{claim}). Therefore, $\ga$ is a duality group
and hence $Mod_S$ is a virtual duality group.

\begin{frem} \label{solomon-tits} {\em
The analogy between the arithmetic and the Teichm\"uller cases breaks down in an important respect.
Recall that the topological dimension of $\mathcal C(S)$ is equal to $3g-4$,
since there are at most $3g-3$ disjoint non-homotopic simple closed curves
in the surface $S$ of genus $g$.
In contrast to the Solomon-Tits Theorem for spherical Tits buildings,
it is not clear why the topological dimension $3g-4$ of $\mathcal C(S)$
is much higher than its homotopy dimension $2g-2$. 
This is quite different from spherical Tits buildings whose topological
and homotopy dimensions  agree.}
\end{frem}

\noindent
{\sc Michigan State University,\\
Department of Mathematics,\\
Wells Hall,\\
East Lansing, MI 48824-1027.\\

\noindent
E-mail:} ivanov@math.msu.edu\\

\noindent
{\sc University of Michigan,\\
Department of Mathematics,\\
East Hall, 530 Church Street,\\
Ann Arbor, MI 48109-1043\\

\noindent
E-mail:} lji@umich.edu

\end{document}